\magnification=1200
\overfullrule=0pt
\centerline {\bf Existence of fixed points for a particular
multifunction}\par
\bigskip
\bigskip
\centerline {BIAGIO RICCERI}\par
\bigskip
\bigskip
\centerline {Department of Mathematics}\par
\centerline {University of Catania}\par
\centerline {Viale A. Doria 6}\par
\centerline {95125 Catania, Italy}\par
\centerline {E-mail: ricceri@dmi.unict.it}
\bigskip
\bigskip
\noindent
{\bf Abstract.} In this paper, we prove that if
$E$ is an infinite-dimensional reflexive real Banach space
possessing the Kadec-Klee property, then, for every compact function
from the unit sphere $S$ of $E$ to the dual $E^*$ satisfying the
condition $\inf_{x\in S}\|f(x)\|_{E^*}>0$, there exists
$\hat x\in S$ such that
$$f(\hat x)(\hat x)=\|f(\hat x)\|_{E^*}\ .$$
\par
\bigskip
\bigskip
\bigskip
\bigskip
Here and in the sequel, $E$ is a reflexive real Banach space, with
dual $E^*$. Set
$$S=\{x\in E : \|x\|=1\}\ .$$
Let $f:S\to E^*$ be a continuous function.\par
\smallskip
In this very short paper, we are interested in the existence of some $\hat x\in S$
such that
$$f(\hat x)(\hat x)=\|f(\hat x)\|_{E^*}\ .$$
Clearly, if, for each $x\in S$,
we set
$$\Phi_f(x):=\{y\in S : f(x)(y)=\|f(x)\|_{E^*}\}\ ,$$
our problem is equivalent to finding a fixed point of the multifunction
$\Phi_f$.
\smallskip
Note that, by reflexivity, we have $\Phi_f(x)\neq
\emptyset$ for all $x\in S$. However, the mere continuity of $f$ is not
enough to guarantee the existence of solutions to our problem.
\smallskip
Let us recall that,
when dim$(E)=\infty$, $E$ is said to have the Kadec-Klee property if
for every sequence $\{x_n\}$ in $S$ weakly converging to $x\in S$,
one has $\lim_{n\to \infty}\|x_n-x\|=0$.\par
\smallskip
Also, we say that a function $\psi:S\to E^*$ is compact if it is
continuous and $\psi(S)$ is relatively compact.\par
\smallskip
Here is our contribution about the above problem.\par
\medskip
THEOREM 1. - {\it Let $E$ be infinite-dimensional and have the
Kadec-Klee property, and
let $f:S\to E^*$ be a compact function such that
$$\inf_{x\in S}\|f(x)\|_{E^*}>0\ . \eqno{(1)}$$
Then, there exists $\hat x\in S$ such that
$$f(\hat x)(\hat x)=\|f(\hat x)\|_{E^*}\ .$$}
\medskip
Let us recall that a multifunction $F:X\to 2^Y$ between two topological spaces
is said to be upper semicontinuous if, for each closed
set $C\subseteq Y$,
the set
$$F^-(C):=\{x\in X : F(x)\cap C\neq\emptyset\}$$
is closed in $X$.\par
\smallskip
{\it Proof of Theorem 1}. Consider the multifunction $\Psi:E^*\to
2^S$ defined by putting
$$\Psi(\varphi)=\{x\in S: \varphi(x)=\|\varphi\|_{E^*}\}$$
for all $\varphi\in E^*$. Let us show that the
restriction of $\Psi$ to $E^*\setminus \{0\}$ is upper semicontinuous.
To this end, let $C\subseteq S$ be a (non-empty)
closed set and let $\{\varphi_n\}$ be
a sequence in $\Psi^-(C)\setminus \{0\}$ converging in $E^*$
to $\varphi\neq 0$. We have to show that $\varphi\in \Psi^-(C)$. For each
$n\in {\bf N}$, choose $x_n\in C$ so that
$$\varphi_n(x_n)=\|\varphi_n\|_{E^*}\ .\eqno{(2)}$$
By reflexivity, there is a subsequence $\{x_{n_k}\}$ weakly converging
to some $x\in E$. Reading $(2)$ with $n_k$ instead of
$n$ and passing to the limit for $k\to \infty$, we get
$$\varphi(x)=\|\varphi\|_{E^*}\ .$$
Since $\varphi\neq 0$, we have $x\in S$. Consequently, since $E$ has the
Kadec-Klee property,
$\{x_{n_k}\}$ converges strongly to $x$. Therefore, since $C$ is closed, we
have $x\in C$. Hence, $x\in \Psi(\varphi)\cap C$ and so $\varphi\in \Psi^-(C)$.
Next, observe that, for each $\varphi\in E^*\setminus \{0\}$, the set
$\Psi(\varphi)$ is bounded, closed and convex, and so it is weakly compact, by
reflexivity. But, since $E$ has the Kadec-Klee property, each weakly compact
subset of $S$ is compact, in view of the Eberlein-Smulyan theorem. So, each
set $\Psi(\varphi)$, with $\varphi\neq 0$, is compact. 
Next, note that,
by $(1)$, $K:=\overline {f(S)}$ is a compact set in $E^*$ which does not contain
$0$. Consequently, by the upper semicontinuity of $\Psi_{|(E^*\setminus \{0\})}$, 
the set $\Psi(K)$ is compact ([2], Theorem 7.4.2). 
Since
dim$(E)=\infty$, there is a continuous function $\omega:B\to S$ such that
$\omega(x)=x$ for all $x\in S$, where $B$ is the closed unit ball of $E$.
Finally, denote by
$Y$ the closed convex hull of $\Psi(K)$ and set
$$G(x)=\Psi(f(\omega(x)))$$
for all $x\in Y$. So, $G$ is an upper semicontinuous multifunction (as
the composition of the upper semicontinuous multifunction
$\Psi_{|(E^*\setminus \{0\})}$ and the continuous function $f\circ \omega$)
with non-empty, closed and convex values, from the compact convex set
$Y$ into itself. Then, by the Fan-Kakutani theorem ([1]), there exists
$\hat x\in Y$ such that $\hat x\in G(\hat x)$. Then, since $\hat x\in S$,
we have $\omega(\hat x)=\hat x$, and so $\hat x\in \Phi_f(\hat x)$, as desired.
\hfill $\bigtriangleup$ \par
\medskip
 Some remarks about the assumptions of Theorem 1 are now in order.\par
\smallskip
Assume that $(E,\langle\cdot,\cdot,\rangle)$ is a Hilbert space. Consider
the continuous function $f:S\to E^*$ defined by putting
$$f(x)(y)=-\langle x,y\rangle$$
for all $x\in S$, $y\in E$.
In this case, we have
$$f(x)(x)=-1$$
and
$$\|f(x)\|_{E^*}=1$$
for all $x\in S$. This example shows, at the same time, that
Theorem 1 is no longer true if either the infinite dimensionality of $E$ or
the compactness of $f$ is removed.\par
\smallskip
Also, note that, concerning the compactness of $f$, a more sophisticated
example can be provided in any infinite-dimensional Banach space.
Actually, in this case, E. Michael ([3]) proved that there
exists a continuous function $f:S\to E^*$ such that
$$f(x)(x)=0$$
and
$$\|f(x)\|_{E^*}=1$$
for all $x\in S$.\par
\smallskip
Concerning the necessity of $(1)$, consider the following example.\par
\smallskip
Let $E$ be the space of all absolutely continuous functions $u:[0,1]\to
{\bf R}$ with $u'\in L^2([0,1])$. In other words, let $E$ be the usual
Sobolev space $H^1(0,1)$, with the usual norm
$$\|u\|=\left ( \int_0^1 |u'(t)|^2dt +\int_0^1 |u(t)|^2dt\right ) ^{1\over 2}\ .$$
Consider the continuous function $f:S\to E^*$ defined by putting
$$f(u)(v)=-\int_0^1 u(t)v(t)dt$$
for all $u\in S$, $v\in E$. Note that $E$ has the Kadec-Klee property since it
is a Hilbert space. Moreover, since $E$ is compactly embedded into $C^0([0,1])$,
the function $f$ is compact. Finally, we have
$$f(u)(u)=-\int_0^1 |u(t)|^2dt<0<\int_0^1 |u(t)|^2dt\leq \|f(u)\|_{E^*}$$
for all $u\in S$.\par
\smallskip
We conclude by proposing the following problem.\par
\medskip
PROBLEM 1. - Let $E$ be an infinite-dimensional reflexive real Banach space
such that, for each compact function $f:S\to E^*$ satisfying $(1)$, there
exists $\hat x\in S$ for which
$$f(\hat x)(\hat x)=\|f(\hat x)\|_{E^*}\ .$$
Then, does $E$ possess the Kadec-Klee property ?
\bigskip
\bigskip
\centerline {\bf References}\par
\bigskip
\bigskip
\noindent
[1]\hskip 5pt K. FAN, {\it Fixed-point and minimax theorems in locally convex topological
linear spaces}, Proc. Nat. Acad. Sci. U.S.A., {\bf 38} (1952), 121-126.\par
\smallskip
\noindent
[2]\hskip 5pt E. KLEIN and A. C. THOMPSON, {\it Theory of correspondences},
John Wiley $\&$ Sons, 1984.\par
\smallskip
\noindent
[3]\hskip 5pt E. MICHAEL, {\it Continuous selections avoiding a set},
Topology Appl., {\bf 28} (1988), 195-213.\par

\bye